\documentclass[preprint,11pt]{elsarticle}
\usepackage[top=2.5cm, bottom=2.5cm, left=2.5cm, right=2cm]{geometry}
\usepackage{amssymb,latexsym,amscd,amsmath,amsfonts,enumerate,supertabular}
\usepackage{graphicx}
\usepackage{tabularx}
\usepackage{subfigure}
\usepackage{enumerate}
 \usepackage{epstopdf}
 \usepackage{color}
 \usepackage{natbib}
 \journal{Elsevier}
\newtheorem{theorem}{Theorem}

\newtheorem{remark}{Remark}

\newenvironment{proof}[1][Proof]{\begin{trivlist}
\item[\hskip \labelsep {\bfseries #1}]}{\end{trivlist}}

\begin{document}
\begin{frontmatter}
\title{{\bf Positive and elementary stable explicit nonstandard Runge-Kutta methods for a class of autonomous dynamical systems
}}
\author[1]{Quang A Dang}
\ead{dangquanga@cic.vast.vn}
\author[2]{Manh Tuan Hoang}
\ead{hmtuan01121990@gmail.com}
\address[1]{Center for Informatics and Computing, Vietnam Academy of Science and Technology (VAST),\\ 
18 Hoang Quoc Viet, Cau Giay, Hanoi, Vietnam}
\address[2]{Institute of Information Technology, 
Vietnam Academy of Science and Technology (VAST), \\
18 Hoang Quoc Viet, Cau Giay, Hanoi, Vietnam}
 \begin{abstract}
\small
In this paper, we construct explicit nonstandard Runge-Kutta (ENRK) methods which have higher accuracy order and preserve two important properties of autonomous dynamical systems, namely, the positivity and linear stability. These methods are based on the classical explicit Runge-Kutta methods, where instead of the usual $h$ in the formulas there stands a function $\varphi (h)$. It is proved that the constructed methods preserve the accuracy order of the original Runge-Kutta methods. 
 The numerical simulations confirm the validity of the obtained theoretical results.
\small
 \end{abstract}
\small 
\small  
\end{frontmatter}
\section {Introduction}
In this paper, we consider the general autonomous system of the following form
\begin{equation}\label{eq:1a}
\dfrac{dy}{dt} = f(y), \qquad y(t_0) = y_0 \in \mathbb{R}^n,
\end{equation}
where $y = \big[y^1, y^2, \ldots , y^n\big]^T:[t_0, T) \to \mathbb{R}^n$, the functions $f = \big[f^1, f^2, \ldots , f^n\big]^T:\mathbb{R}^n \to \mathbb{R}^n$ is differentiable. A one-step numerical scheme with a step size $h$, that approximates the solution $y(t_k)$ of System \eqref{eq:1a} can be written as $D_h(y_k) = F_h(f; y_k)$, where $D_h(y_k) \approx {dy}/{dt}, F(f; y_k) \approx f(y)$, and $t_k = t_0 + kh$.\par
Our first objective is to construct difference schemes preserving the linear stability of the equilibrium points 
of  System \eqref{eq:1a} for all finite step-size $h > 0$. These schemes are called also 
 elementary stable \cite{AL, DK1, DK2}. It should be emphasized that 
 standard finite difference schemes cannot preserve properties of the differential equations for any step-sizes 
  $h > 0$, including the linear stability. Mickens called this phenomenon numerical instability
 \cite{Mickens2}.\par
 The construction of elementary stable difference schemes play especially important role in 
 numerical solution of differential equations and numerical simulation of nonlinear dynamical systems. In 2005, Dimitrov and Kojouharov \cite{DK1} proposed a method for constructing elementary stable NSFD methods for general two-dimensional autonomous dynamical systems. These NSFD methods are based on the  explicit and implicit Euler  and the second-order Runge-Kutta methods. Later, in 2007 these results are extented for the general n-dimensional dynamical systems, namely, NSFD schemes preserving elementary stability are constructed based on the $\theta$-methods and the second-order Runge-Kutta methods \cite{DK2}. One important action in the construction of the elementary stable NSFD schemes   is the replacement of the  standard denominator function $\varphi(h) = h$ by the nonstandard denominator function, which is bounded from above for any  $h > 0$. This nonstandard denominator function guarantees the difference schemes to be elementary stable but it influences the order of accuracy of original Runge-Kutta methods. In general, the order of accuracy of the obtained nonstandard NSFD methods is lower than the order of accuracy of the original Runge-Kutta methods. Usually, NSFD schemes, which are dynamically consistent with differential equations, have first order accuracy \cite{Mickens2}. This motivates the problem of increasing the order of accuracy of NSFD schemes. Recently, some higher order  NSFD schemes were constructed \cite{Chen-Charpentier1, MVaquero1, MVaquero}. They are based on the combination of nonstandard schemes and Richardson's extrapolation. Besides, it is worthy to mention the interesting design of NSFD schemes of second order accuracy for differential equations with polynomial right-hand sides based on the combination of the exact schemes for some special polynomials \cite{Chen-Charpentier2}.  \par
Differently from the above ways, in this paper we construct higher order nonstandard finite methods based on Runge-Kutta methods without extrapolation. More specifically, we design a class of elementary stable higher order NSFD methods based on explicit standard Runge-Kutta (ESRK) methods with some correction of the usual denominator $h$. Namely, the usual denominator $h$ is replaced by a nonstandard denominator function $\varphi(h)$ which is bounded from above by a number $\varphi^*$, i.e., $\varphi(h) < \varphi^*$ for all $h > 0$. The number $\varphi^*$ plays a role as a threshold for elementary stability of
 ESRK method and it is explicitly determined for each method with particular problems .\par
Many processes and phenomena in  fields of application are described by dynamical systems which possess additional important properties including the positivity. The construction of difference schemes preserving the both positivity and linear stability  is very important but not simple. There are several works concerning this topic, among them the typical works are of Dimitrov and Kojouharov \cite{DK3, DK4}, Wood and Kojouharov \cite{Wood}.\par
{For many dynamical systems with relevant right-hand sides, it is possible to determine 
 a positivity step size thresholds of Runge-Kutta methods (see \cite{Horvath1, Horvath2}), that is a number $H > 0$ which guarantees that the numerical solutions obtained by Runge-Kutta methods are positive for any step size  $h < H$. Then, it is clear that the replacement of the denominator functions not only guarantees the elementary stability but also guarantees the positivity of Runge-Kutta methods. Namely, we need only to choose the function $\varphi(h)$ satisfying the conditions
   $\varphi(h) < \tau^* : = \min\{\varphi^*, H\}$, for all $h > 0$, where $\varphi^*$ and $H$ 
   are the elementary stability threshold and the positivity threshold of ESRK, respectively.
   This number $\tau^*$ (if it is determined) is called the positivity and elementary stable step size threshold, or briefly PES step size thresholds of  ESRK methods.}\par
In this paper, first a PES step size threshold for ESRK methods are determined on the base of the stability analysis of ESRK methods and the application of results for positive Runge-Kutta methods \cite{Horvath1, Horvath2}. After that the choice of nonstandard denominator functions for preserving the accuracy order of the original ESRK methods will be investigated. In result we shall obtain ENRK methods of higher order of accuracy, which preserve simultaneously the positivity and elementary stability of dynamical systems. For short these methods will called positive and elementary stable (PESN) methods (see \cite{DK3, DK4}). Our result resolves the contradiction between the dynamical consistency with higher order of accuracy of NSFD schemes. Here,
  we do not use implicit Runge-Kutta methods because they require computational cost for solving a nonlinear system of algebraic equations for determining the stage-coefficients $K_i$ at each iteration. Moreover, in order to guarantee the existence of a unique solution of the nonlinear system the step size $h$ must be not large \cite[Theorem 7.2]{Hairer}. Besides, the determination of positivity and elementary stable step size thresholds of explicit methods, in general, is  simpler.\par
The  paper is organized as follows. In Section 2, we construct positive and elementary stable ENRK methods. The influence of the denominator function will be analyzed in Section 3. Section 4 reports numerical simulations and some applications of the constructed ENRK methods. The last Section is discussion and some conclusions.
\section{Positive and Elementary stable ENRK methods}
\subsection{Elementary stable ENRK methods}
In this Subsection we construct elementary stable NSFD schemes based on the s-stage ESRK methods \cite{Ascher, Hairer}. For this purpose we consider NSFD schemes of the form
\begin{equation}\label{eq:1}
\begin{split}
&K_1 = f(y_k),  \quad K_2 = f(y_k + \varphi(h) a_{21}K_1), \quad \ldots \quad  K_s = f(y_k + \varphi(h) \sum_{j = 1}^{s - 1}a_{sj}K_j ),\\
&\dfrac{y_{k + 1} - y_k}{\varphi(h)} = b_1K_1 + b_2K_2 + \ldots + b_sK_s, \qquad \qquad \varphi(h) = h + \mathcal{O}{(h^2)}, \quad h \to 0.
\end{split}
\end{equation}
For simplicity in presentation we omit the argument 
$h$ in the function $\varphi(h)$ and use the notation $(A, b^T, \varphi)$ to refer to ENRK methods \eqref{eq:1}, where $A = (a_{ij})_{s \times s}$ and $b = (b_1, \ldots, b_s)$. Notice that $(A, b^T, h)$ is common ESRK methods.\par
Assume that system \eqref{eq:1a} has a finite number of equilibria and $Re (\lambda) \ne 0$ for $\lambda \in \Omega$, where $\Omega = \bigcup_{y^* \in \Gamma}\sigma(J(y^*))$, and $\Gamma$ represents the set of all equilibria of system \eqref{eq:1a}. We shall consider two cases of order of accuracy $p$ and number of stages $s$ of the original  ESRK  $(A, b^T, h)$. They are the case $p = s$ and the case $p < s$. 
Notice that the case $p= s$ holds only for $p \leq 4$  because according to \cite[Section 4.3]{Ascher}, \cite[Chapter II]{Hairer} the order of s-stage ESRK methods with $s \ge 5$ always is less than $s$.
\begin{theorem}\label{maintheorem}
Suppose that the original ESRK  $(A, b^T, h)$ has the order of accuracy $p = s$. Then there exists a number  $\varphi^* := \varphi^*(p, \Omega) > 0$, playing the role of the elementary stability threshold of 
 \eqref{eq:1}, i.e. \eqref{eq:1} is elementary stable if $\varphi(h)$ satisfies the condition
\begin{equation}\label{eq:2}
0 < \varphi(h) < \varphi^*, \quad \forall h > 0.
\end{equation}
\end{theorem}
\begin{proof}
Suppose $y^*$ is an equilibrium point of  \eqref{eq:1a}. Denote by $J = J(y^*)$ and $\hat{J} = \hat{J}(y^*)$  the Jacobians of \eqref{eq:1a} and \eqref{eq:1}, respectively. Based on the proofs in \cite{DK1, DK2} and the results of the stability function of ESRK \cite[Section 4.4]{Ascher}, \cite[Section IV.2]{Hairer1} we can deduce that if $\lambda_i$ $(i = \overline{1, n})$ are the eigenvalues of  $J$ then the respective eigenvalues 
 $\mu_i$ $(i = \overline{1, n})$ of $\hat{J}$  are
\begin{equation}\label{eq:3}
\mu_i = \sum_{j = 0}^p\dfrac{z^j}{j!}= \dfrac{z^p}{p!} + \dfrac{z^{p-1}}{(p-1)!} + \ldots + \dfrac{z^2}{2} + z + 1, \qquad i = \overline{1, n}, \qquad z = \varphi \lambda_i.
\end{equation}
In other words, the method
 $(A, b^T, \varphi)$ maps the eigenvalues $\lambda_i$ to the eigenvalues $\mu_i$, respectively. Notice that the right hand side of \eqref{eq:3} is the stability function of an 
   ESRK method having the order of accuracy $p = s$.\\
Suppose that  $\lambda = Re(\lambda) + i Im(\lambda)\in \sigma(J)$ is a complex eigenvalue represented in the trigonometric form
\begin{equation*}
\lambda = r_{\lambda}\Big(\cos(\theta_{\lambda}) + i \sin(\theta_{\lambda})\Big), \, r_{\lambda} = |\lambda|, \, \cos(\theta_{\lambda}) = {Re(\lambda)}/{|\lambda|}, \, \sin(\theta_{\lambda}) = {Im(\lambda)}/{|\lambda|},
\end{equation*}
then in accordance with $\lambda$ we have
 \begin{equation*}
\mu_\lambda = \dfrac{(\varphi \lambda)^p}{p!} + \dfrac{(\varphi \lambda)^{p-1}}{(p-1)!} + \ldots + \varphi \lambda + 1 = 
\sum_{j = 0}^{p}\dfrac{1}{j!}\varphi^j r_{\lambda}^j \cos(j \theta_{\lambda}) + i \sum_{j = 1}^{p}\dfrac{1}{j!}\varphi^j r_\lambda^j \sin(j \theta_{\lambda}),
\end{equation*}
\begin{equation*}
\begin{split}
|\mu_\lambda|^2 &= \Big(\sum_{j = 0}^{p}\dfrac{1}{j!}\varphi^j r_{\lambda}^j \cos(j \theta_{\lambda}) \Big)^2 + \Big(\sum_{j = 1}^{p}\dfrac{1}{j!}\varphi^j r_{\lambda}^j \sin(j \theta_{\lambda})\Big)^2\\
& =  \Big(\sum_{j = 1}^{p}\dfrac{1}{j!}\varphi^j r_{\lambda}^j \cos(j \theta_{\lambda}) \Big)^2 + \Big(\sum_{j = 1}^{p}\dfrac{1}{j!}\varphi^j r_{\lambda}^j \sin(j \theta_{\lambda})\Big)^2
+ 2 \sum_{j = 2}^{p}\dfrac{1}{j!}\varphi^j r_{\lambda}^j \cos(j \theta_{\lambda}) + 2 \varphi r_{\lambda} \cos(\theta_{\lambda}) +  1.
\end{split}
\end{equation*}
From here it follows that $|\mu_\lambda|^2 < 1$ if and only if $\mathcal{P}_{2p - 1}(\varphi, \lambda) < 0$, where
\begin{equation}\label{eq:5}
\begin{split}
\mathcal{P}_{2p - 1}(\varphi, \lambda) := \Big(\sum_{j = 1}^{p}\dfrac{1}{j!}{(\varphi r_{\lambda})}^{j - 1/2} \cos(j \theta_{\lambda}) \Big)^2 &+ \Big(\sum_{j = 1}^{p}\dfrac{1}{j!}{(\varphi r_{\lambda})}^{j - 1/2} \sin(j \theta_{\lambda})\Big)^2\\
&+ 2 \sum_{j = 2}^{p}\dfrac{1}{j!}{(\varphi r_{\lambda})}^{j - 1}\cos(j \theta_{\lambda}) + 2\cos(\theta_{\lambda}) < 0.
\end{split}
\end{equation}
Obviously, $\mathcal{P}_{2p - 1}(\varphi, \lambda)$ is a polynomial of  $\varphi$ having the highest degree $2p-1$ with coefficients depending on $r_\lambda$ and $\theta_\lambda$, and the free term is $2\cos(\theta_\lambda)$. \par 
Consider two cases of stability of   $y^*$:\\
{\bf Case 1.} $y^*$ is a linearly stable equilibrium point of \eqref{eq:1a}. Then $Re(\lambda_i) < 0$ for any $\lambda_i \in \sigma(J)$ $(i = \overline{1, n})$. The necessary and sufficient condition for $y^*$ to be a linearly stable equilibrium point of \eqref{eq:1} is $|\mu_i| < 1$ for any $\mu_i \in \sigma(\hat{J})$ $(i = \overline{1, n})$. Since $Re(\lambda) < 0$ there holds $\cos(\theta_{\lambda}) < 0$. On the other hand, due to $\lim_{\varphi \to 0} \mathcal{P}_{2p - 1}(\varphi, \lambda) = 2\cos(\theta_{\lambda}) < 0$ from the definition of limit of a function it follows that there exists a number 
  $\varphi_1 > 0$ such that $\mathcal{P}_{2p - 1}(\varphi, \lambda) < 0$ for all $\varphi \in (0, \varphi_1)$, where $\mathcal{P}_{2p - 1}(\varphi, \lambda)$ is defined by \eqref{eq:5}. \\  
{\bf Case 2.} $y^*$ is a linearly unstable equilibrium of \eqref{eq:1a}. Then there exists $ \lambda_i \in \sigma(J)$ for some $1 \leq i \leq n$ such that $Re(\lambda_i) > 0$. The necessary and sufficient condition for $y^*$ to be a linearly unstable equilibrium point of \eqref{eq:1} is the existence of some $j$ such that $|\mu_j| > 1$. Suppose for $i = l$ there holds $Re(\lambda_l) > 0$.  The necessary and sufficient condition for $|\mu_l| > 1$ is $\mathcal{P}_{2p - 1}(\varphi, \lambda_l) > 0$, where $\mathcal{P}_{2p - 1}(\varphi, \lambda)$ is defined by \eqref{eq:5}. Since $Re(\lambda_l) > 0$ we have $\cos(\theta_{\lambda_l}) > 0$. Therefore, there exists a number $\varphi_2 > 0$ such that $\mathcal{P}_{2p - 1}(\varphi, \lambda_l) > 0$ for all $\varphi \in (0, \varphi_2)$.\par
Denote $\Omega^+ = \bigcup_{y^* \in \Gamma^+}\sigma(J(y^*))$, $\Omega_- = \Big\{\xi \in \bigcup_{y^* \in \Gamma_-}\sigma(J(y^*)): Re(\xi) > 0 \Big\}$, where $\Gamma^+$ and $\Gamma_-$ are the set of linearly stable equilibria and the set of linearly unstable equilibria of 
 \eqref{eq:1a}, respectively. 
 {Set ${\varphi^*} = \min\Big\{\varphi^+, \varphi_-\Big\}$, where}
\begin{equation*}\label{eq:e1}
\begin{split}
\varphi^+(\lambda) &= \sup_{\varphi^+ >  0}\bigg\{\varphi^{+}: \mathcal{P}_{2p - 1}(\varphi, \lambda) < 0, \, \forall \varphi \in (0, \varphi^+), \, \lambda \in \Omega^+ \bigg\}, \quad \varphi^+ = \min\bigg\{\varphi^+(\lambda): \lambda \in \Omega^+ \bigg\},\\
\varphi_-(\lambda) &= \sup_{\varphi_{-} >  0}\bigg\{\varphi_-: \mathcal{P}_{2p - 1}(\varphi, \lambda) > 0, \forall \varphi \in (0, \varphi_-),\, \lambda \in \Omega_-\bigg\}, \quad  \varphi_- = \min\bigg\{\varphi_-(\lambda): \lambda \in \Omega_- \bigg\}.
\end{split}
 \end{equation*}
Then from  Case 1 and Case 2 it follows that if $0 < \varphi < \varphi^*$ then \eqref{eq:1} is elementary stable. The proof of the theorem is complete.
\end{proof}
The following theorem for the case $p < s$ is stated and proved in a similar way as Theorem \ref{maintheorem}.
\begin{theorem}\label{maintheorem2}
Suppose that the original ESRK method $(A, b^T, h)$ has the order of accuracy $p < s$. Then there exists a number $\varphi^* := \varphi^*(p, \Omega, A, b^T, s) > 0$ such that \eqref{eq:1} is elementary stable if the function $\varphi(h)$ satisfies the condition $0 < \varphi(h) < \varphi^*$ for all  $h > 0$.
\end{theorem}
\begin{proof}
Since $p < s$ the eigenvalues $\mu_i$ $\in \hat{J}$ $(i = \overline{1, n})$ corresponding to the the eigenvalues $\lambda_i \in J$ are defined by
\begin{equation}\label{eq:3a}
\mu_i = \sum_{j= 0}^p\dfrac{z^j}{j!} + \sum_{j > p}^s z^j b^T A^{j - 1} \textbf{1} = 1 + z + \ldots + \dfrac{z^p}{p!} + \sum_{j > p}^s a_{j}z^j, \quad z = \varphi \lambda_i, \quad a_j := a_j(p, A, b^T, s).
\end{equation}
Notice that the right hand side of  \eqref{eq:3a} is the stability function of an ESRK method having the order of accuracy $p <s$ 
{ \cite[Section 4.4]{Ascher}, \cite[Section IV.2]{Hairer1}}. Repeating the proof of Theorem \ref{maintheorem} with the attention that now $\mathcal{P}_{2p - 1}(\varphi, \lambda)$ defined by \eqref{eq:5} is replaced by the polynomial $\mathcal{P}_{2s - 1}(\varphi, \lambda)$ of the degree $2s- 1$ we come to the conclusion of the theorem.
\end{proof}
\begin{remark}
It is well-known that for ensuring linear stability ESRK methods have to use step-sizes small enough, i.e., ESRK methods are conditionally stable. Meanwhile, according to Theorems \ref{maintheorem} and \ref{maintheorem2} ENRK methods with the appropriately chosen denominator function are unconditionally stable. Notice that if $p = s$ then $\varphi^*$ depends only on $p$ but if $p <s$ then it depends also on $A, b^T, s$. From the proof of the theorems it is possible to design an algorithm for determining the number $\varphi ^*$ based on the finding of the minimal positive root of the polynomials $\mathcal{P}_{2p - 1}(\varphi, \lambda)$ or $\mathcal{P}_{2s - 1}(\varphi, \lambda)$. In the cases if  $p = s = 1$ and $p = s =  2$, the explicit formulas for $\varphi^*$ are given in \cite{DK1, DK2}.
\end{remark}
\subsection{Positive ENRK methods}
In this Subsection, we propose a method for constructing positive ENRK methods based on the results of positivity of Runge-Kutta methods \cite{Horvath1, Horvath2}.\par
Suppose that the right-hand side $f(y)$ of the equation \eqref{eq:1a} satisfies conditions such that the solution of the equation \eqref{eq:1a} is nonnegative for all $y_0 \geq 0$. The set of such functions $f$  is denoted by $\mathcal{P}$ (see \cite{Horvath1, Horvath2}). The necessary and sufficient conditions  for $f \in \mathcal{P}$ can be found in \cite{Horvath1, Horvath2}. It is possible to find  a positivity step size thresholds of Runge-Kutta methods having $R(A, b) > 0$ in many cases of  $f$, for example, if $f$ belongs to one of the sets $\mathcal{F}^*$, $\mathcal{F}^*(\rho)$, $\mathcal{F}_\infty^*(\rho)$, $\mathcal{F}_{\infty, w}^*(\rho)$ or $\mathcal{P}_\alpha$ \cite{Horvath2}. Here $R(A, b)$ is positivity radius of Runge-Kutta methods with the coefficient scheme $(A, b)$ \cite{GERISCH, Horvath1, Horvath2, KRAAIJEVANGER}. The radius $R(A, b)$ is used by Kraaijevanger \cite{KRAAIJEVANGER} in the study of contractivity of RK methods and also used in the nonlinear positivity theory for RK methods by Horvath \cite{Horvath1, Horvath2}. The results related to the properties of  $R(A, b)$ may be found in \cite{KRAAIJEVANGER}.\par
Now suppose that the right-hand side $f$ belongs to one of the sets $\mathcal{F}^*$, $\mathcal{F}^*(\rho)$, $\mathcal{F}_\infty^*(\rho)$, $\mathcal{F}_{\infty, w}^*(\rho)$ or $\mathcal{P}_\alpha$.
Notice that there are many important differential equations models have right-hand sides belonging to the above sets, for example,  metapopulation models \cite{Amarasekare, Keymer}, predator-prey systems \cite{DK5}, computer virus model with graded cure rates \cite{Yang}, vaccination model with multiple endemic states \cite{KV}, mathematical biology system \cite{Allen}\ldots Then for explicit Runge-Kutta methods with $R(A, b) > 0$, we can determine a positivity step size thresholds depending on $R(A, b)$ (see \cite[Proposition 3 and Theorem 4]{Horvath2}). Suppose this number is $H > 0$. Combining this fact with Theorem \ref{maintheorem} and Theorem \ref{maintheorem2} we have ESRK methods \eqref{eq:1} which are PES method for \eqref{eq:1a} if
\begin{equation}\label{eq:8a}
\varphi(h) < \tau^* := \min\{\varphi^*, H\}, \, \forall \, h > 0.
\end{equation}
Notice that the above conditions for existence of a positivity step size thresholds $H > 0$ is only sufficient conditions for the positivity of Runge-Kutta methods. For many Runge-Kutta methods, although the number $H > 0$ cannot be determined by this way, there exists a number
$H^* > 0$ such that the method is positive for any  $h < H^*$. Besides, the number  $H > 0$ determined by the mentioned method may be not strict positivity step size thresholds (see \cite{Horvath2}). In general, from the estimate $|y(t_k) - y_k| = \mathcal{O}(h^p) \, as \, h \to 0$, it is possible to conjecture that when  $h$ is sufficiently small, that is, $h < H^*$ for some $H^*$,  Runge-Kutta methods are positive. However, to determine the threshold 
 $H^*$ is difficult.\par
Next, based on the results of the number $R(A, b)$ \cite[Section 9]{KRAAIJEVANGER}, we choose some explicit Runge-Kutta methods having $R(A, b) > 0$ for determining positivity step size thresholds.
\begin{enumerate}
\item For $1$-stage methods only the Euler method has $R(A, b) = 1$.
\item For $2$-stages method  we choose the second order Heun method (explicit trapezoidal method - RK2) having $R(A, b) = 1$.
\item Among $4$-stages methods, we choose the 4-stages of order 3 (RK43) having maximal $R(A, b)$.  This is the method defined by \cite[Section 9]{KRAAIJEVANGER} with $R(A, b) = 2$:
\begin{equation*}
b_1 = b_2 = b_3 = 1/6, \qquad b_4 = 1/2, \qquad a_{21} = 1/2, \qquad a_{31} = a_{32} = 1/2, \qquad a_{41} = a_{42} = a_{43} = 1/6.
\end{equation*}
\item According to \cite[Theorem 9.6]{KRAAIJEVANGER}, there does not exist  explicit 4-stage coefficient scheme $(A, b)$ with classical order $p = 4$ and $R(A, b)$ > 0.  Therefore, we choose 5-stage method with $p = 4$ (RK54). This method has maximal $R(A, b) \approx 1. 50818$ (see \cite[p. 522, Section 9]{KRAAIJEVANGER}).
\end{enumerate}
\begin{remark}
The condition $R(A, b)> 0$ makes narrow possible explicit Runge-Kutta methods. For example, two well-known 4-stages methods, namely, 
classical Runge-Kutta method and ${3}/{8}$-rule do not belong to the methods having
$R(A, b) > 0$. However, in numerical simulations it will be seen that these sufficient conditions may be freed. For Runge-Kutta methods, which not necessarily satisfy these conditions, it is possible to determine the positivity threshold.
\end{remark}
\section{The influence of the denominator function}
In this Section, we make an analysis of the influence of denominator function with the aim to show  some choices of this function for preserving the accuracy order of the original  ESRK.\par
First, obviously, the standard denominator function $\varphi(h) = h$ does not satisfy \eqref{eq:8a}. It is easy to choose the function $\varphi(h)$ satisfying \eqref{eq:8a}, for example (see \cite{AL, DK1, DK2, Mickens2})
\begin{equation}\label{eq:e1}
\varphi_1(h) = \dfrac{1 - e^{-\tau_1 h}}{\tau_1}, \qquad\tau_1 > 0.
\end{equation}
However, the change of denominator function may cause the decrease of the accuracy order of ENRK, i.e., it does not preserve the accuracy order of the original  ESRK. Therefore, it is important to choose the function $\varphi (h)$ so that the accuracy order of the original  ESRK is preserved.
\begin{theorem}
Suppose the original ESRK  $(A, b^T, h)$ has accuracy order $p$. Then ENRK methods $(A, b^T, \varphi)$ also have accuracy order $p$ if $\varphi(h)$ satisfies the condition
\begin{equation}\label{eq:6}
\varphi(h) = h + \mathcal{O}(h^{p + 1}), \qquad h \to 0.
\end{equation}
\end{theorem}
\begin{proof}
First, notice that if $\varphi(h)$ satisfies \eqref{eq:6} then $\varphi'(0) = 1$ and $\varphi(0) = \varphi''(0) = \varphi^{(3)}(0) = \ldots = \varphi^{(p)}(0) = 0$. Applying the method for construction of the system of order conditions for ESRK based on Taylor expansion \cite[Chapter II]{Hairer}, it is easy to deduce that the system of conditions for order $p$ for $(A, b^T, h)$ and $(A, b^T, \varphi(h))$ coincide if $\varphi(h)$ satisfies \eqref{eq:6}. Thus, the theorem is proved.
\end{proof}
Now, to choose the function $\varphi(h)$ satisfying simultaneously \eqref{eq:8a} and \eqref{eq:6} we consider the class of functions 
\begin{equation}\label{eq:e2}
\varphi_2(h) = h e^{-\tau_2 h^m}, \qquad m \in \mathbb{Z}_+, \qquad m \geq p, \qquad \tau_2 > 0.
\end{equation}
Clearly, $\varphi_2(h)$ satisfies \eqref{eq:6} and reaches the maximal value at $h^* = \sqrt[m]{{1}/{m\tau_2}}$, i.e., for any $h > 0$ we have $\varphi_2(h) \leq \varphi(h^*) =  e^{-{1}/{m}}{\sqrt [m]{{1}/{m\tau_2}}} \to 0$ as $\tau_2 \to \infty$. It means that there always exists $\tau_2 > 0$ such that $\varphi_2(h)$ satisfies \eqref{eq:8a}. Here it suffices to choose $\tau_2 > \tau_{opt}^2 := [{m e (\tau^*)^m}]^{-1}$. From \eqref{eq:e1} it follows that the condition for  $\varphi_1(h)$ to satisfy \eqref{eq:8a} is $\tau_1 > \tau_{opt}^1 := {(\tau^*)}^{-1}$.\par
It is easy to see that $\varphi_i(\tau_i) \to h$ as $\tau_i \to 0$ $(i = 1, 2)$. It means that when $h$  and $\tau_i$ are small then ENRK methods (with the denominator functions $\varphi_i(h)$) have the accuracy order equal to that of the original ESRK. In other words, in order to ensure the accuracy order  $\tau_i$ must be chosen as small as possible. However, the choice of $\tau_i$ depends on the value of $\tau_{opt}^i$.\par
It is seen that the constraint $\tau_1 > \tau_{opt}^1 := {(\tau^*)}^{-1}$ does not allow to choose $\tau_1$ arbitrarily small, especially when $\tau^*$ is very small. It is the reason leading to the decrease of the accuracy order of the original ESRK when choosing the function $\varphi_1(h)$ defined by \eqref{eq:e1}. The numerical simulations in Section 4 will clearly demonstrate this fact.\\
Concerning the function $\varphi_2(h)$ we see that if $\tau^* \geq 1$ then $\tau_{opt}^2 = [{m e (\tau^*)^m}]^{-1}  \to 0$ as $m \to \infty$, therefore, if $m$ is chosen large enough then $\tau_{opt}^2$ is very small. When $\tau^* < 1$ then $\tau_{opt}^2  \to \infty$ as $m \to \infty$, hence the choice $m = p$ is best. If $\tau^*$ is very small (e.g. for stiff problems \cite{Ascher, Hairer1}) then $\tau_{opt}^2$ is also very large. Then ESRK methods should to use the step-size $h < \tau^*$ for guaranteeing the stability and the accuracy order. Since $\varphi_2(h)$ satisfies \eqref{eq:6} ENRK methods also have the same accuracy order as the original  ESRK when $h < \tau^*$. When $h \geq \tau^*$  ESRK methods become unstable while ENRK methods stay stable. Thus, if $h$ is small then  $\varphi_2(h)$ gives the accuracy better than $\varphi_1(h)$.\par
However, since $\varphi_2(h) \to 0$ as $h \to \infty$, whenever $h$ is large $\varphi_2(h) \approx 0$. Then the numerical solution obtained at all steps are slightly different from the initial value. Therefore, it is impossible to estimate the asymptotic behavior of the solution. In other words, for large $h$ it is not recommended to use the function $\varphi_2(h)$.  In this case, the function  $\varphi_1(h)$ is advantageous because it is a monotonically increasing and upper bounded function. 
As $h \to \infty$ there holds $\varphi_1(h) \to 1/\tau_1 < \tau^*$, therefore, although $h$ is large enough the numerical solution remains to have linear stability as the exact solution (see \cite{DK1, DK2}).\par
 In order to overcome the shortcoming of $\varphi_2(h)$, we propose to use a class of new functions $\varphi(h)$, which guarantee the accuracy order when $h$ is small and guarantee the asymptotic behavior of numerical solution when $h$ is large. It is the class of functions of the form
\begin{equation}\label{eq:10}
\varphi_3(h) = \theta(h)\varphi_2(h) + \big(1 - \theta(h)\big)\varphi_1(h),
\end{equation}
where the function $\theta(h)$ has the property $\theta(h) = 1 + \mathcal{O}(h^p)$ as $h \to 0$; $0 < \theta(h) < 1$ for all $h > 0$; $\lim_{h \to 0} \theta(h)= 1$ and $\lim_{h \to \infty}\theta(h) = 0$.\par
\begin{figure}[!ht]
\centering
\includegraphics[height=8cm,width=8cm]{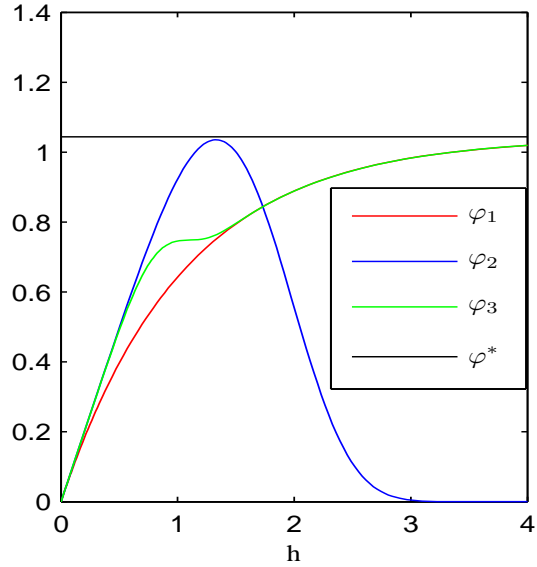}
 \caption{Graphs of the functions $\varphi_i(h)$.}\label{fig1}
\end{figure}

Obviously, $\varphi_3(h)$ satisfies \eqref{eq:6} and if $\varphi_1(h)$ and $\varphi_2(h)$ satisfy \eqref{eq:8a} then so does $\varphi_3(h)$. Especially, when $h$ is small then $\varphi_3(h)$ is equivalent to $\varphi_2(h)$, so  this guarantees the best error,  conversely, when $h$ is large then $\varphi_3(h)$ is equivalent to $\varphi_1(h)$ and this guarantees the asymptotic behavior of numerical solution.
It is best to choose   $\theta(h) = e^{-h^k}$, where $k \in \mathbb{Z}_+$. Figure \ref{fig1} depicts the graphs of the function $\varphi_i(h) (i = \overline{1, 3})$ for some particular values of the parameter.
\section{Numerical simulations and some applications of the ENRK methods}
In this Section we apply the constructed RK methods to some important applied models. Numerical experiments confirm the validity of the obtained theoretical results of ENRK methods.\par
\subsection{A predator-prey system with Beddington-DeAngelis functional response}
We consider the following predator-prey system with Beddington-DeAngelis functional response \cite{DK5} which was considered in \cite{DK2}
\begin{equation}\label{eq:example1}
\dfrac{dx}{dt} = x - \dfrac{Axy}{1 + x + y}, \qquad \dfrac{dy}{dt} = \dfrac{Exy}{1 + x + y} - Dy, \qquad (x(0), y(0)) = (1, 1.6),
\end{equation}
where $x$ and $y$ represent the prey and predator population sizes, respectively, and the values of the constants are $A = 2, D = 1$ and $E = 10$. Mathematical analysis of System \eqref{eq:example1} shows that there exist two equilibria $(0, 0)$ and $(\dfrac{AD}{AE - E - AD}, \dfrac{E}{AE - E - AD}) = (0.25, 1.25)$, among them the equilibrium $(0.25, 1.25)$ is globally stable in the interior of the first quadrant and the equilibrium $(0, 0)$ is unstable. The eigenvalues of $J(0, 0)$ are given by $\lambda_1 = 1$ and $\lambda_2 = -1$, and the eigenvalues of $J(0.25, 1.25)$ are given by 
$\lambda_3 = -\dfrac{1}{5} + \dfrac{3}{5}i = r\Big(\cos(\theta) + i\sin(\theta)\Big)$ and $\lambda_4 =  -\dfrac{1}{5} - \dfrac{3}{5}i = r\Big(\cos(-\theta) + i\sin(-\theta)\Big)$, where $r = \sqrt{\dfrac{2}{5}}$ and $\theta = 0.6024\pi$.\par
Since it is impossible to find the exact solution of the system we use the numerical solution obtained by a $11$-stage Runge-Kutta method of accuracy order 8 (RK8)
 \cite{Cooper} with  $h = 10^{-5}$ as a benchmark solution. Here $error = max_{k}\big\{|x_k - X_k| + |y_k - Y_k|\big\}$ is used as a measure for the accuracy of   ENRK methods, where $\{(x_k, y_k)\}$ and $\{(X_k, Y_k)\}$ are the solutions obtained by ENRK methods and the benchmark solution, respectively. Besides, 
 $rate := log_{h_1/h_2}(error(h1)/error(h2))$ (see \cite[Example 4.1]{Ascher}) is an approximation for accuracy order of the methods.  \par
In this example we shall consider 
ENRK1 (based on the Euler method), ENRK2 (based on the second Heun method), ENRK43 (based on the RK43 method), ENRK54 (based on the RK54 method). The numbers $R(A, b)$ for these methods are  1, 1, 2, and 1.50818, respectively (see Subsection 2.2). Besides,
 ENRK4 and (based on the classical 4-stages RK method) are also considered although these methods do not possess $R(A, b) > 0$.\\
It is easy to verify that the right-hand side $f$ of the model belongs to the set $\mathcal{P}_\alpha : = \{f| f(t, v) + \alpha v \geq 0 \, \, \text{for \,all }\, \, t, v \geq 0\}$ with $\alpha = \max\{A -1 , D\} = 1$ (see \cite{Horvath1, Horvath2}), therefore it is easy to determine
positivity step size thresholds for ENRK methods.\par
  Besides, it is not difficult to determine the polynomials
 $\mathcal{P}_{2p -1}$ and $\mathcal{P}_{2s -1}(\varphi, \lambda_i) \; (i = \overline{1, 4})$ corresponding to ENRK methods. From them it is easy to determine the numbers
   $\varphi^*$. The numbers $\varphi^*$, $\tau_{opt}^i$ $(i=1, 2)$ and the functions $\varphi_j(h)$ $(j = \overline{1, 3})$ for the methods ENRK1, ENRK2, ENRK43, ENRK54, ENRK4  are given in Table \ref{tabl1}. The errors and the rates of the methods for small $h$ are reported in Tables \ref{tabl2}-\ref{tabl6}, where for short, the columns with the headings $\varphi (h)$ and $\varphi_i$ stand for the errors of the methods with the denominator function $\varphi (h)$ and $\varphi_i$, respectively; the columns with the headings $rate_i$ stand for the rates of the methods corresponding to $\varphi_i$.\\
Notice that it is impossible to determine the positive threshold for  ENRK4 because $R(A, b) = 0$. In this case it is only possible to determine the elementary stability threshold for ENRK4. However, many numerical simulations show that the elementary stability threshold is also the positive threshold. This fact is also seen from the numerical simulations in \cite{DK2}.
\begin{table}
\centering
\setlength{\tabcolsep}{0.12cm}
\caption{The values $\tau_{opt}^i$ and the denominator functions $\varphi_i(h)$ $(i = \overline{1, 3})$ of the ENRK methods.}\label{tabl1}
\medbreak
\begin{tabular}{ c c  c  c c c c c c c}
\hline
Methods&$(s, p)$&$\varphi^*$&$H$&$\tau^*$&$\tau^1_{opt}$&$\varphi_1(h)$&$\tau^2_{opt}$&$\varphi_2(h)$&$\theta(h)$ in \eqref{eq:10}\\
\hline
\\
ENRK1&(1, 1)&0.9998&1&0.9998&1.0002&$\dfrac{1- e^{-1.0005h}}{1.0005}$&0.0919&$he^{-0.095h^4}$&$e^{-0.01h^2}$\\
\\
ENRK2&(2, 2)&2.6604&1&1&1&$1- e^{-h}$&0.0920&$he^{-0.095h^4}$&$e^{-0.01h^4}$\\
\\
ENRK43&(4, 3)&4.7332&2&2&0.5&$\dfrac{1 - e^{-0.55}}{0.55}$&9.5802e-004&$he^{-0.001h^6}$&$e^{-h^6}$\\
\\
ENRK54&(5, 4)&5.0631&1.50818&1.50818&0.6631&$\dfrac{1 - e^{-0.68}}{0.68}$&1.7179e-003&$he^{-0.002h^8}$&$e^{-h^8}$\\
\\
ENRK4&(4, 4)&4.4476&*&4.4476&0.2248&$\dfrac{1 - e^{-0.25}}{0.25}$&7.9214e-006&$he^{-0.0001h^6}$&$e^{-0.01}h^6$\\
\\
\hline
\end{tabular} 
\end{table}
\begin{table}
\centering
\setlength{\tabcolsep}{0.12cm}
\caption{The errors and rates of ENRK1 methods.}\label{tabl2}
\medbreak
\begin{tabular}{ c c c c c c c c}
\hline
$h$&$\varphi(h) = h$&$\varphi_1$&$rate_1$&$\varphi_2$&$rate_2$&$\varphi_3$&$rate_3$\\
\hline
$0.2$&0.4303&0.6056&& 0.4304&&0.4304&\\
$0.1$&0.2032&0.2937&1.0443&0.2032 &1.0827&0.2032 &1.0827\\
$0.05$&0.0986 &0.1444&1.0239&0.0986&1.0439&0.0986&1.0439 \\
$0.01$&0.0192&0.0285&1.0091&0.0192&1.0156&0.0192& 1.0156 \\
$0.005$&0.0096&0.0142&1.0027&0.0096&1.0045&0.0096&1.0045\\
$0.001$&0.0019&0.0028&1.0009&0.0019&1.0016&0.0019&1.0016\\
\hline
\end{tabular} 
\end{table}
\begin{table}
\centering
\setlength{\tabcolsep}{0.12cm}
\caption{The errors and rates of ENRK2 methods.}\label{tabl3}
\medbreak
\begin{tabular}{ c c c c c c c c}
\hline
$h$&$\varphi(h) = h$&$\varphi_1$&$rate_1$&$\varphi_2$&$rate_2$&$\varphi_3$&$rate_3$\\
\hline
$0.2$& 7.3223e-003  &4.1755e-001&&7.1013e-003& & 7.0992e-003&\\
$0.1$&1.7189e-003& 2.1136e-001& 0.9823& 1.7052e-003& 2.0581&1.7051e-003& 2.0578 \\
$0.05$& 4.1773e-004&1.0622e-001 &  0.9926 &4.1687e-004&2.0323& 4.1686e-004&   2.0322\\
$0.01$&1.6354e-005&  2.1321e-002&  0.9978&1.6352e-005& 2.0121& 1.6352e-005&  2.0121 \\
$0.005$&4.0770e-006& 1.0665e-002& 0.9994& 4.0769e-006& 2.0040&4.0769e-006& 2.0040\\
$0.001$&1.6271e-007&2.1337e-003& 0.9998&1.6271e-007&  2.0014&1.6271e-007& 2.0014 \\
\hline
\end{tabular} 
\end{table}
\begin{table}
\centering
\setlength{\tabcolsep}{0.12cm}
\caption{The errors and rates of ENRK43 methods.}\label{tabl4}
\medbreak
\begin{tabular}{ c c c c c c c c}
\hline
$h$&$\varphi(h) = h$&$\varphi_1$&$rate_1$&$\varphi_2$&$rate_2$&$\varphi_3$&$rate_3$\\
\hline
$0.2$& 5.8286e-004 & 1.9063e-001&& 5.8275e-004& & 5.7796e-004&\\
$0.1$&7.1911e-005& 9.5672e-002&0.9946&7.1910e-005& 3.0186&7.1872e-005&  3.0075 \\
$0.05$& 8.9428e-006& 4.7924e-002 &  0.9973 &8.9428e-006& 3.0074& 8.9425e-006&   3.0067\\
$0.01$&7.1300e-008&  9.5989e-003& 0.9991&7.1300e-008& 3.0021& 7.1300e-008&  3.0021 \\
$0.005$&8.9081e-009& 4.8003e-003& 0.9997& 8.9081e-009& 3.0007&8.9081e-009&3.0007\\
$0.001$&7.1181e-011&9.6021e-004& 0.9999&7.1181e-011& 3.0007&7.1181e-011&  3.0007 \\
\hline
\end{tabular} 
\end{table}

\begin{table}
\centering
\setlength{\tabcolsep}{0.12cm}
\caption{The errors and rates of ENRK54 methods.}\label{tabl5}
\medbreak
\begin{tabular}{ c c c c c c c c}
\hline
$h$&$\varphi(h) = h$&$\varphi_1$&$rate_1$&$\varphi_2$&$rate_2$&$\varphi_3$&$rate_3$\\
\hline
$0.2$& 3.1359e-005 & 2.8632e-001&& 3.1368e-005& &  3.1665e-005&\\
$0.1$&2.0695e-006& 1.4419e-001&0.9897&2.0695e-006 & 3.9219&2.0700e-006&  3.9352 \\
$0.05$&1.3274e-007& 7.2338e-002 &  0.9952 &1.3274e-007&  3.9626& 1.3274e-007&   3.9629\\
$0.01$&2.1686e-010& 1.4502e-002& 0.9985& 2.1686e-010&  3.9871& 2.1686e-010&  3.9871 \\
$0.005$&1.3706e-011&7.2531e-003&  0.9996& 1.3706e-011& 3.9839&1.3706e-011&3.9839\\
$0.001$& 1.9159e-012& 1.4510e-003& 0.9999&1.9159e-012&  1.2225& 1.9159e-012& 1.2225 \\
\hline
\end{tabular} 
\end{table}

\begin{table}
\centering
\setlength{\tabcolsep}{0.12cm}
\caption{The errors and rates of ENRK4 methods.}\label{tabl6}
\medbreak
\begin{tabular}{ c c c c c c c c}
\hline
$h$&$\varphi(h) = h$&$\varphi_1$&$rate_1$&$\varphi_2$&$rate_2$&$\varphi_3$&$rate_3$\\
\hline
$0.2$&1.9481e-005& 1.0622e-001&& 1.9488e-005&&2.1385e-005&\\
$0.1$& 1.1945e-006 & 5.3233e-002 &0.9967& 1.1946e-006 &  4.0280& 1.2044e-006&  4.1502\\
$0.05$&7.3021e-008&2.6646e-002&0.9984&7.3022e-008& 4.0321&7.3099e-008&   4.0423\\
$0.01$&1.1429e-010 & 5.3336e-003& 0.9995& 1.1429e-010&  4.0137&1.1430e-010&  4.0143 \\
$0.005$&7.2312e-012 &2.6671e-003&0.9998&7.2312e-012&3.9824& 7.2312e-012 & 3.9824\\
$0.001$&1.9159e-012&5.3346e-004&0.9999&1.9159e-012&  0.8253&1.9159e-012& 0.8253\\
\hline
\end{tabular} 
\end{table}
From Tables \ref{tabl2}-\ref{tabl6} we see that the replacement  of the denominator function
 $\varphi(h) = h$ by the function $\varphi_1(h)$ decreases the accuracy order of the original ESRK.
 Namely, ENRK methods have only order $1$. Meanwhile, ENRK methods with the appropriate denominator functions $\varphi_i(h)$ ($i = 2, 3)$ have the accuracy orders equal to those of the original ESRK methods. In other words, the accuracy orders of the original 
 ESRK methods are preserved. 
 In the columns $rate_2$ and $rate_3$ of Tables \ref{tabl5} and \ref{tabl6} we see an unexpected phenomenon, when $h$ is small the rates decrease. A similar phenomenon also was indicated in \cite[Example 4.1]{Ascher} when studying explicit standard Runge-Kutta methods.  
 The reason of this is that  the rounding errors generally increase as $h$ decreases. \par
\begin{figure}[!ht]
\centering
\includegraphics[height=10cm,width=12cm]{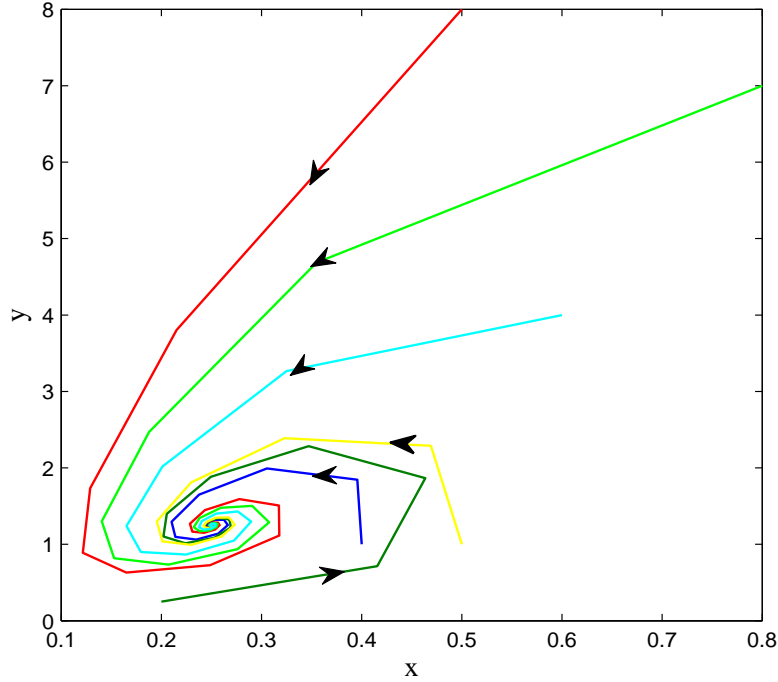}
\caption{Numerical solution obtained by  ENRK54 with $\varphi_3(h)$ and h = 4.}\label{fig2}
\end{figure}

{Next}, we consider ENRK54 methods for large $h$, specifically  $h = 4$. Then $\varphi_2(4) \approx 4.7667e-057$, therefore it is possible to accept $\varphi_2(4) = 0$.  
The computational experiments show that the numerical solution obtained when using the function $\varphi_2(h)$ at all grid nodes are constant and are equal to the initial values. Meanwhile, the numerical solutions obtained by the methods with the use of  the functions $\varphi_1(h)$ and $\varphi_3(h)$ have the asymptotic behavior similar to that of the exact solution
 Figure \ref{fig2} depicts the numerical solution obtained by 
  ENRK54 with $\varphi_3(h) = e^{-h^8}he^{-0.002h^6} + (1 - e^{-h^8})\dfrac{1 - e^{-0.68h}}{0.68}$ and $h = 4$. 
 The experiments for the other methods give similar results. This completely agrees with the analysis made in  Section 3. Notice that from
Figure \ref{fig2} we see that the global stability of the model is also preserved. \par
Finally, for the purpose of comparison, we apply the method proposed by Wood and Kojouharov \cite{Wood} for the model. This is a NSFD scheme preserving the positivity and elementary stability of the dynamical system based on nonlocal approximation.
For this method the positivity step size threshold is $H = \infty$, and its the elementary stability threshold is the same of Euler method.
Therefore, following the results in \cite{Wood} we choose the denominator function $\varphi(h) = \dfrac{1- e^{-1.0005h}}{1.0005}$. The error and convergence rate of the method are reported in Table \ref{tabl7}. From the table we see that the method has accuracy order 1 equal the order of the explicit Euler method but the Euler method is somewhat better (see Table \ref{tabl2}). 
\begin{table}
\centering
\setlength{\tabcolsep}{0.12cm}
\caption{The errors and rates of the Wood and Kojouharov methods.}\label{tabl7}
\medbreak
\begin{tabular}{ c c c c c c c c}
\hline
$h$&$\varphi(h) = h$&$\varphi_1$&$rate_1$&$\varphi_2$&$rate_2$&$\varphi_3$&$rate_3$\\
\hline
$0.2$& 0.5390& 0.8143&& 0.5391&&0.5621 &\\
$0.1$& 0.2575  & 0.4041 &1.0107& 0.2576&  1.0658&0.2918&  0.9460\\
$0.05$&0.1257&0.2008&1.0093&0.1257 & 1.0348& 0.1538 &   0.9237\\
$0.01$&0.0247& 0.0399& 1.0038& 0.0247& 1.0121&0.0326&  0.9644\\
$0.005$&0.0123&0.0199&1.0012&0.0123&1.0036& 0.0164&0.9889\\
$0.001$&0.0025&0.0040&1.0004&0.0025&  1.0012&0.0033& 0.9962\\
\hline
\end{tabular} 
\end{table}
\subsection{A vaccination model with multiple endemic states}
We consider the following vaccination model with multiple endemic states \cite{KV} which was considered in \cite{DK2}
\begin{equation*}
\begin{split}
\dfrac{dS}{dt} &= \mu N - \beta SI/N - (\mu + \phi)S + cI + \delta V,\\
\dfrac{dI}{dt} &= \beta SI/N - (\mu + c)I,\\
\dfrac{dV}{dt}&=\phi S - (\mu + \delta)V,
\end{split}
\end{equation*}
where the constants $\beta = 0.7$, $c = 0.1$, $\mu = 0.8$, $\delta = 0.8$ and $\phi = 0.8$. In the above
model the total (constant) population size $N = 100$ is divided into three classes susceptibles ($S$), infectives ($I$) and vaccinated ($V$) and it is assumed that the vaccine is completely effective in preventing infection. Mathematical analysis of this System  shows that the disease free equilibrium $(S^*, I^*, V^*) = \Big(\dfrac{(\mu + \delta)N}{\mu + \delta + \phi}, 0, \dfrac{\phi N}{\mu + \delta + \phi}\Big) = \Big(\dfrac{200}{3}, 0, \dfrac{100}{3}\Big)$ is globally asymptotically stable \cite{KV}. The eigenvalues of  $J(S^*, I^*, V^*)$  are given by $\lambda_1 = -0.8$, $\lambda_2 = -2.4$, $\lambda_3 = -\dfrac{13}{30}$. Therefore, it is easy to determine the elementary stability threshold for ENRK methods. Moreover, the right-hand side of the system belongs to the set 
 $\mathcal{P}_\alpha$ where $\alpha = \max\{\beta + \mu + \phi, \mu + c, \mu + \delta \} = 2.5$, hence we can determine the positive thresholds for ENRK. Table \ref{tabl8} gives the positive and elementary stability thresholds for ENRK.

\begin{figure}[!ht]
\centering
\includegraphics[height=11cm,width=12cm]{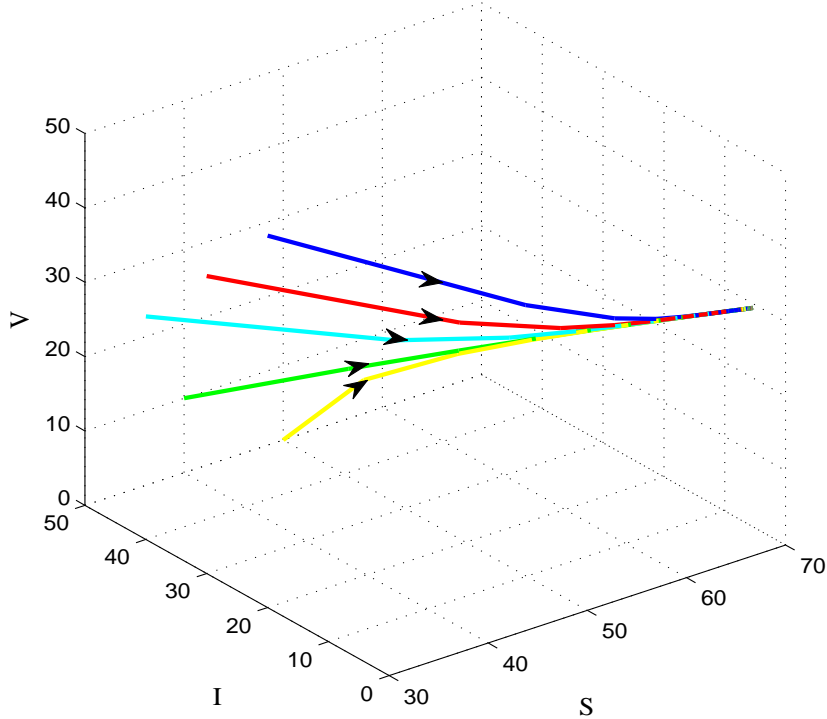}
\caption{Numerical solution obtained by ENRK54 với $\varphi_3(h) = e^{-h^6}he^{-0.5h^4} + ( 1-e^{-h^6})\dfrac{1 - e^{-1.6h}}{1.6} $ and $h = 2$.}\label{fig3}
\end{figure}

\begin{table}
\centering
\setlength{\tabcolsep}{0.12cm}
\caption{ Positive and elementary stability thresholds for ENRK.}\label{tabl8}
\medbreak
\begin{tabular}{ c c c c  c c}
\hline
Methods&ENRK1&ENRK2&ENRK43&ENRK54&ENRK4\\
\hline
$\varphi^*$& 0.8333&0.8333&2.1499& 2.2068&1.1605\\
$H = R(A, b)/\alpha$&0.4&0.4&0.8&0.6631&*\\
$\tau^* = \min\{\varphi^*, H\}$&0.4&0.4&0.8&0.6631&*\\
\hline
\end{tabular} 
\end{table}
Based on these results we choose the denominator function for ENRK methods. The numerical solution obtained by ENRK54 is depicted in Figure \ref{fig3}. We see that the global stability of the model is preserved, while the numerical simulations 
 in \cite{DK2} show that the second order RK2 and Euler method do not preserve this property of the model. In addition, the advantage of the ENRK methods is that they preserve  the property of the model for all $h > 0$ and have higher order of accuarcy  when $h$ is small.
\subsection{Metapopulation models}
Consider the metapopulation model proposed by Keymer in 2000  \cite{Keymer}. It is described by the system of three nonlinear differential equations
\begin{equation}\label{eq:meta1}
\dfrac{dp_0}{dt} = e(p_1 + p_2) - \lambda p_0, \quad \dfrac{dp_1}{dt} = \lambda p_0 - \beta p_1 p_2 + \delta p_2 - ep_1, \quad \dfrac{dp_2}{dt}= \beta p_1 p_2 - (\delta + e)p_2, 
\end{equation}
where $p_0$, $p_1$ and $p_2$ denote the proportion of uninhabitable patches, the proportion of the habitable patches that are not occupied and the proportion of habitable patches that are occupied, respectively, $\lambda$ is the rate of patch creation, $e$ is the rate of patch destruction, $\delta$ is the rate of population extinction and $\beta$ is the rate of propagule reproduction.\par
Recently, we constructed a discrete model, which is dynamically consistent with the model
\eqref{eq:meta1} based on NSFD schemes \cite{DQA1, DQA2}. The related numerical simulations show that RK2, RK4, Euler methods do not preserve important properties of the model including the positivity and linear stability of equilibrium points.\par
It is easy to see that the right-hand side of \eqref{eq:meta1} belongs to $\mathcal{P}_\alpha$, where $\alpha = \max\{\lambda, \beta + e, \delta + e\}$. Therefore, in a similar way, we can construct  ENRK methods preserving the positivity and stability of equilibrium points. Moreover, the numerical simulations in \cite{DQA1, DQA2} hint that if the positivity and stability of the model are preserved then other important properties, including the invariant property $p_0(t) + p_1(t) + p_2(t) \equiv 1$ , are also preserved by ENRK.\par
In a similar way, consider the metapopulation model proposed by Amarasekare and Possingham in 2001 \cite{Amarasekare}. It is described by the system of four nonlinear differential equations
\begin{equation}\label{eq:meta2}
\begin{split}
&\dfrac{dI}{dt} = \beta_I SI - e_I I + fL - gI, \qquad \qquad \dfrac{dS}{dt} = e_I I - \beta_I SI + fR - gS,\\
&\dfrac{dL}{dt} = gI - fL - e_LL + \beta_LRI, \qquad \qquad \dfrac{dR}{dt} = gS - fR + e_LL - \beta_LRI.\\
\end{split}
\end{equation}
Here $f$ is the disturbance frequency and $g$, the rate of habitat succession. Quantities $e_I$ and $e_L$ represent local extinction rates, and $\beta_I$ and $\beta_L$ the per patch colonization rates of infected and latent patches, respectively. The total number of patches in the system is assumed to be constant such that $I + S + L + R = P$. Alternatively, $I, S, L$ and $R$ can be thought of as the frequency of each patch type in the landscape in which case $I + S + L + R = 1$.\par
The model \eqref{eq:meta2} has analogous properties as the model \eqref{eq:meta1}. Recently, we also constructed NSFD preserving important properties of the model. It is worthy to notice that the stability properties of the model were established with the help of an extention of Lyapunov theorem \cite{DQA2}.\\
It is easy to see that the right-hand side of \eqref{eq:meta2} also belongs to $\mathcal{P}_\alpha$, where $\alpha = \max\{e_I + g, \beta_I + g, f + e_L, f + \beta_L\}$. Therefore, we can construct ENRK preserving the positivity and stability properties of the model. These methods not only preserve important properties of the model but also have higher order of accuracy.
\section{Discussion and Conclusions}
In this paper, starting from explicit standard Runge-Kutta methods we have constructed explicit nonstandard Runge-Kutta (ENRK) methods, which are PES methods and preserve the accuracy order of the original Runge-Kutta methods. When the step-size is small these ENRK methods have higher order of accuracy and when the step-size is large they preserve the stability of the dynamical system. The numerical simulations confirm the validity of the obtained theoretical results of the constructed methods. 
 \par 
 The key in the construction of elementary stable ENRK methods is the determination of the positive and elementary stability thresholds of the methods. After the determination of these thresholds it is followed by the replacement of
the standard denominator by suitable nonstandard denominator functions (Section 3). The numerical simulations in Section 4 show that not only positivity and stability properties but also other important properties of the model including the global stability are preserved. \par
Theorem 1 and Theorem 2 may be considered as an extension of the previous results of 
Dimitrov and Kojouharov \cite{DK1, DK2}. Meanwhile, the positivity threshold is determined based on the results of positivity of RK methods \cite{Horvath1, Horvath2}. The conditions for the positivity threshold by this way narrows RK methods designed for considered dynamical systems although the performed numerical simulations show that these conditions may be freed. Therefore, it is needed to wide the conditions for the positivity threshold. 
Moreover, the analysis made in Section 3 hints that it is reasonably to choose the positivity and elementary stability threshold  be possibly maximal (strict). Meanwhile, the elementary stability threshold determined by Theorem 1 and 2 is  strict, the positivity threshold determined via   $R(A, b)$ may be not strict. Therefore, the problem of choosing maximal threshold is to be investigated.\par
Some researches agree that when the denominator function is bounded from above by some appropriate  value  $\varphi^*$ then NSFD schemes are dynamically consistent with continuous models (see, e.g. \cite{DK3, DK4, DQA1, DQA2, Chen-Charpentier1} \ldots). This value $\varphi^*$ depends on the continuous systems, their properties and on each NSFD scheme. For example, in this paper, for ensuring the elementary stability, this values   $\varphi^*$ is determined explicitly with the help of  Theorem \ref{maintheorem} and Theorem \ref{maintheorem2}.\par
In general, from the convergence order $p$ of ESRK methods ($|y(t_k) - y_k| = \mathcal{O}(h^p) \, as \, h \to 0$), it is predicted that for sufficiently small  $h$  ESRK methods preserve some properties of the continuous models (for example, the positivity and stability of equilibrium points). It means that there exists $h^* > 0$ such that for any $h < h^*$  ESRK schemes preserve the important properties of continuous models. Then the choice of the denominator functions satisfying  $\varphi(h) < h^*$ for all $h > 0$ ensure that ENRK schemes are dynamically consistent with continuous models. Therefore, it is possible to extend the results of this paper to construct ENRK dynamically consistent with the differential equations   models. Besides, these results are applicable to some systems of ordinary differential equations obtained after discretization of evolution PDEs by the method of lines. \par
In the future we will develop the proposed technique to construct NSFD methods which have higher accuracy order and are dynamically consistent with autonomous dynamical systems. In parallel, the application and extension of these methods for relevant problems also will be investigated.
\section*{Acknowledgments}
This work is supported by Vietnam National Foundation for Science and Technology Development (NAFOSTED) under the grant  number 102.01-2017.306.

\end{document}